\documentclass[a4paper, 10pt]{article}
\usepackage{amssymb}

\usepackage{amsmath,amsthm}
\usepackage{multicol}

\begin{document}

\title{Swarming Behavior of Multi-Agent Systems %
\thanks{%
Peking University Swarm Research Group Members: Professor Long
Wang, Professor Tianguang Chu, Ms. Shumei Mu, Miss Hong Shi and
Miss Bo Liu. Corresponding author: Professor Long Wang.}}
\author{P. K. U. Swarm \\
Intelligent Control Laboratory, Center for Systems and Control,\\
Department of Mechanics and Engineering Science, Peking University,\\
Beijing 100871, P. R. China, Email: longwang@mech.pku.edu.cn}
\date{}
\maketitle

\begin{minipage}{5.9in}

\textbf{Abstract}: In this paper we consider a continuous-time
anisotropic swarm model in $n$-dimensional space with an
attraction/repulsion function and study its aggregation
properties. It is shown that the swarm members will aggregate and
eventually form a cohesive cluster of finite size around the swarm
center. Moreover, the numerical simulations show that all agents
will eventually enter into and remain in a bounded region around
the swarm center. The model is more general than isotropic swarms
and our results provide further insight into the effect of the
interaction pattern on individual motion in a swarm system.


\end{minipage}
\begin{multicols}{2}

\section{Introduction}
In nature swarming can be found in many organisms ranging from
simple bacteria to more advanced mammals. Examples of swarms
include flocks of birds, schools of fish, herds of animals, and
colonies of bacteria. Such collective behavior has certain
advantages such as avoiding predators and increasing the chance of
finding food. Recently, there has been a growing interest in
biomimicry of forging and swarming for using in engineering
applications such as optimization, robotics, military applications
and autonomous air vehicle \cite{K. Passino}--\cite{M. Pachter and
P. Chandler}. Modeling and exploring the collective dynamics has
become an important issue and many papers have appeared \cite{Y.
Liu 1}--\cite{A. Czirok 2}. However, results on the anisotropic
swarms are relatively few. The study of anisotropic swarm is very
difficult though the anisotropic swarm is a ubiquitous phenomenon,
including natural phenomena and social phenomena.

Gazi and Passino \cite{V. Gazi and K. M. Passino} proposed an
isotropic swarm model and studied its aggregation, cohesion and
stability properties. Subsequently, Chu and Wang \cite{T. Chu}
generalized their model, considering an anisotropic swarm model,
and obtained the properties of aggregation, cohesion and
completely stability. The coupling matrix $W$ considered in
\cite{T. Chu} is symmetric, that is, the interactions between two
individuals are reciprocal. In this paper, we try to study the
behavior of anisotropic swarms when the coupling matrix is
completely nonsymmetric. The model and the results given here
extend the work on isotropic swarms \cite{V. Gazi and K. M.
Passino} and anisotropic swarms \cite{T. Chu} to more general
cases and further illustrate the effect of the interaction pattern
on individual motion in swarm system.

In the next section we specify an ``individual-based"
continuous-time anisotropic swarm model in an $n$-dimensional
Euclidean space which includes the isotropic model of \cite{V.
Gazi and K. M. Passino} as a special case. Then, under some
assumptions, we show that the swarm can exhibit aggregation in
Section 3. In Section 4, we extend the results in Section 3,
considering a more general attraction/repulsion function. In
Section 5, under some assumptions, we provide some numerical
simulations of the agent motion. We briefly summarize the results
of the paper in Section 6.

\section{Anisotropic Swarms}

We consider a swarm of $N$ individuals (members) in an
$n$-dimensional Euclidean space. We model the individuals as
points and ignore their dimensions. We consider the equation of
motion of individual $i$ described by
\begin{equation}
{\dot{x}^i}=\sum_{j=1}^Nw_{ij}f(x^i-x^j), \;i=1,\cdots,N,
\label{eq1}
\end{equation}
where $x^i\in R^n$ represents the position of individual $i$;
$W=[w_{ij}]\in R^{N\times N}$ with $w_{ij}\geq 0$ for all
$i,j=1,\cdots,N$ is the coupling matrix; $f(\cdot)$ represents the
function of attraction and repulsion between the members. In other
words, the direction and magnitude of motion of each member is
determined as a weighted sum of the attraction and repulsion of
all the other members on this member. The attraction/repulsion
function that we consider is
\begin{equation} f(y)=-y\bigg(a-b\exp\bigg(-\frac{\|y\|^2}{c}\bigg)\bigg), \label{eq2}
\end{equation}
where $a, b$, and $c$ are positive constants such that $b>a$ and
$\|y\|$ is the Euclidean norm given by $\|y\|=\sqrt{y^Ty}$.

In the following discussion we always assume $w_{ii}=0, i=1,
\cdots, N$ in model (1). Moreover, we assume that there are no
isolated clusters in the swarm, that is, $W+W^T$ is irreducible.

Note that the function $f(\cdot)$ is the social potential function
that governs the interindividual interactions and is attractive
for large distances and repulsive for small distances. By equating
$f(y)=0$, one can find that $f(\cdot)$ switches sign at the set of
points defined as ${\cal{Y}}=\big\{y=0$ or
$\|y\|=\delta=\sqrt{c\ln{(b/a)}}\big\}$. The distance $\delta$ is
the distance at which the attraction and repulsion balance. Such a
distance in biological swarms exists \cite{K. Warburton and J.
Lazarus}. Note that it is natural as well as reasonable to require
that any two different swarm members could not occupy the same
position at the same time.

\textbf{Remark 1}: The anisotropic swarm model given here includes
the isotropic model of \cite{V. Gazi and K. M. Passino} as a
special case. Obviously, the present model (1) is more close to
actuality and more meaningful.

\section{Main Results}

In this section, the main results concerning aggregation and
cohesiveness of the swarm (1) are presented. In fact, it is
interesting to investigate collective behavior of the system
rather than to ascertain detailed behavior of each individual. And
due to complex interactions among the multi-agents, in general, it
is very difficult or even impossible to study the specific
behavior of each agent.

Define the center of the swarm members as
$\overline{x}=\frac{1}{N}\sum_{i=1}^Nx^i$, then we have
$\dot{\overline{x}}=\frac{1}{N}\sum_{i=1}^N\dot{x}^i.$ If the
coupling matrix $W$ is symmetric, by the symmetry of $f(\cdot)$
with respect to the origin, the center $\overline{x}$ is
stationary for all $t$ \cite{T. Chu} and the swarm described by
Eqs. (1) and (2) is not drifting on average. Note, however, that
the swarm members may still have relative motions with respect to
the center while the center itself stays stationary and the
members will move toward the swarm center and form a cohesive
cluster around it. However, if the coupling matrix $W$ is
nonsymmetric, the center $\overline{x}$ may not be stationary. An
interesting issue is whether the members will form a cohesive
cluster and which point they will move around. We will deal with
this issue in the following theorem.

\textbf{Theorem 1}: Consider the swarm described by the model in
(1) with an attraction/replusion function $f(\cdot)$ as given in
(2). Assume for any agent $i$, we have
$\sum_{j=1}^Nw_{ij}=\sum_{j=1}^Nw_{ji}$. Then, all agents will
eventually enter into and remain in the bounded region
\begin{equation}
\Omega=\bigg\{x:
\sum_{i=1}^N\|x^i-\overline{x}\|^2\leq{\rho^2}\bigg\}, \label{eq3}
\end{equation}
where$$\rho=\frac{2bM\sqrt {2c}\exp(-\frac{1}{2})}{a\lambda_2};$$
and $\lambda_2$  denotes the second smallest real eigenvalue of
the matrix $L+L^T$; and $M=\sum\limits_{i,j=1}^Nw_{ij}$;
$L=[l_{ij}]$ with
\begin{equation}
\begin{array}{lcl}
&  & {l_{ij}}=\left\{
\begin{array}{l}
-w_{ij}, \\
\sum_{k=1,k\neq{i}}^Nw_{ik},%
\end{array}
\begin{array}{l}
\;i\neq j, \\
\;i=j;
\end{array}
\right.  \\
\end{array}
\label{eq4}
\end{equation}
$\Omega$ provides a bound on the maximum ultimate swarm size.

\begin{proof}
Let $e^i=x^i-\overline{x}$. By the definition of the center
$\overline{x}$ of the swarm and the assumption of
$\sum_{j=1}^Nw_{ij}=\sum_{j=1}^Nw_{ji}$, we have
$$
\dot{\overline{x}}
=\frac{b}{N}\sum\limits_{i=1}^N\bigg[\sum\limits_{j=1}^Nw_{ij}(x^i-x^j)\exp\bigg(-\frac{\|x^i-x^j\|^2}{c}\bigg)\bigg].
$$
To estimate $e^i$, we let $V=\sum_{i=1}^NV_i$ be the Lyapunov
function for the swarm, where $V_i=\frac{1}{2}e^{iT}e^i$.
Evaluating its time derivative along solution of the system (1),
we have
\begin{equation*}
\begin{array}{rl}
\dot{V} =&
\!\!\!-a\displaystyle\sum\limits_{i=1}^N\sum\limits_{j=1}^Nw_{ij}e^{iT}(e^i-e^j)\\
+ & \!\!\! b\displaystyle\sum\limits_{i=1}^Ne^{iT}\bigg\{\sum\limits_{j=1}^Nw_{ij}\beta_{ij}(x^i-x^j)\\
- & \!\!\!\displaystyle\frac{1}{N}\sum\limits_{i=1}^N\bigg[\sum\limits_{j=1}^Nw_{ij}\beta_{ij}(x^i-x^j)\bigg]\bigg\}\\
\leq & \!\!\!-ae^T(L\otimes I)e\\
+ & \!\!\!b\displaystyle\sum\limits_{i=1}^N\sum\limits_{j=1}^Nw_{ij}\beta_{ij}\|x^i-x^j\|\|e^i\|\\
+ & \!\!\!\displaystyle\frac{b}{N}\sum\limits_{i=1}^N\bigg[\sum\limits_{i=1}^N\sum\limits_{j=1}^Nw_{ij}\beta_{ij}\|x^i-x^j\|\|e^i\|\bigg],\\
\end{array}
\end{equation*}%
where $e=(e^{1T},\cdots,e^{NT})^T$ and
$\beta_{ij}=\exp\big(-\frac{\|x^i-x^j\|^2}{c}\big)$, $L\otimes I$
is the Kronecker product of $L$ and $I$ with $L$ as defined in Eq.
(4) and $I$ the identity matrix of order $n$.

Note that each of the functions
$\exp\big(-\frac{\|x^i-x^j\|^2}{c}\big)\|x^i-x^j\|$ is a bounded
function whose maximum occurs at $\|x^i-x^j\|=\sqrt{c/2}$ and is
given by $\sqrt{c/2}\exp(-(1/2))$. Substituting this in the above
inequality and using the fact that $\|e^i\|\leq \sqrt{2V}$, we
obtain
\begin{equation}
\dot{V}\leq -ae^T(L\otimes
I)e+2bM\sqrt{c}\exp\big(-\frac{1}{2}\big)V^{\frac{1}{2}}.
\label{eq5}
\end{equation}

To get further estimate of $\dot{V}$, we only need to estimate the
term $e^T(L\otimes I)e$. Since
$$e^T(L\otimes I)e=\frac{1}{2}e^T\big((L+L^T)\otimes I\big)e,$$ we should analyze
$e^T((L+L^T)\otimes I)e$.  First considering the matrix $L+L^T$
and $L$ as defined in Eq. (4), we have
$L+L^T=[\widetilde{l}_{ij}]$, where
\begin{equation}
\begin{array}{lcl}
&  & {\widetilde{l}_{ij}}=\left\{
\begin{array}{l}
-w_{ij}-w_{ji}, \\
2\sum_{k=1,k\neq{i}}^Nw_{ik},%
\end{array}
\begin{array}{l}
\;i\neq j, \\
\;i=j.
\end{array}
\right.  \\
\end{array}
\label{eq6}
\end{equation}

Using the conditions $\sum_{j=1}^Nw_{ij}=\sum_{j=1}^Nw_{ji}$, we
can conclude that $\lambda=0$ is an eigenvalue of $L+L^T$ and
$u=(l,\cdots,l)^T$ with $l\neq 0$ is the associated eigenvector.
Moreover, since $L+L^T$ is symmetric and $W+W^T (L+L^T)$ is
irreducible, it follows from matrix theory \cite{R. Horn and C. R.
Johnson} that $\lambda=0$ is a simple eigenvalue and all the rest
eigenvalues of $L+L^T$ are real and positive . Therefore, we can
order the eigenvalues of $L+L^T$ as $0=\lambda_1< \lambda_2\leq
\lambda_3\leq \cdots \leq \lambda_n$. Also it is known that the
identity matrix $I$ has an $n$ multiple eigenvalues $\mu=1$ and
$n$ independent eigenvectors
\begin{equation*}
u^1=\left[
\begin{array}{cc}
1 \\
0\\
\vdots\\
0%
\end{array}%
\right] ,\;\;u^2=\left[
\begin{array}{c}
0 \\
1\\
\vdots\\
0%
\end{array}%
\right] ,\;\;\cdots,\;\; u^n=\left[
\begin{array}{cc}
0\\
0\\
\vdots\\
1%
\end{array}%
\right] .
\end{equation*}%

By matrix theory \cite{R. Horn and C. R. Johnson}, the eigenvalues
of $(L+L^T)\otimes I$ are $\lambda_i\mu=\lambda_i$ ($n$ multiple
for each $i$). Next, we consider the matrix $(L+L^T)\otimes I$.
$\lambda=0$ is an $n$ multiple eigenvalues and the associated
eigenvectors are
$$v^1=[u^{1T}, \cdots, u^{1T}]^T, \cdots, v^n=[u^{nT}, \cdots, u^{nT}]^T.$$

Therefore , $e^T\big((L+L^T)\otimes I\big)e=0$ implies that $e$
must lie in the eigenspace of $(L+L^T)\otimes I$ spanned by
eigenvectors $v^1, \cdots, v^n$ corresponding to the zero
engenvalue , that is, $e^1=e^2=\cdots=e^N$. This occurs only when
$e^1=e^2=\cdots=e^N=0$, but this is impossible for the swarm
system under consideration, because it implies that the $N$
individuals occupy the same position at the same time. Hence, for
any solution $x$ of system (1), $e$ must be in the subspace
spanned by eigenvectors of $(L+L^T)\otimes I$ corresponding to the
nonzero eigenvalues. Then, $e^T\big((L+L^T)\otimes I\big)e\geq
\lambda_2\|e\|^2=2\lambda_2 V$. From (5), we have
\begin{equation*}
\begin{array}{rl}
\dot{V}\leq& \!\!\!-a\lambda_2V
+2bM\sqrt{c}\exp(-\frac{1}{2})V^{\frac{1}{2}}\\
= & \!\!\!-\bigg[a\lambda_2V^{1/2}-2bM\sqrt{c}\exp(-\frac{1}{2})\bigg]V^{\frac{1}{2}}\\
< & \!\!\!0%
\end{array}
\end{equation*}%
whenever
$$V(x)>\bigg(\frac{2bM\sqrt{c}\exp(-1/2)}{a\lambda_2}\bigg)^2.$$
Therefore, any solution of system (1) will eventually enter into
and remain in $\Omega$.
\end{proof}

Theorem 1 shows that the swarm members will aggregate and form a
bounded cluster around the swarm center.

\textbf{Remark 2}: The above discussions explicitly show the
effect of the coupling matrix $W$ on aggregation and cohesion of
the swarm.

\textbf{Remark 3}: The conditions given in the above theorem
$\sum_{j=1}^Nw_{ij}=\sum_{j=1}^Nw_{ji}$ include the case as a
special case when the coupling matrix $W$ is a symmetric matrix.

\textbf{Remark 4}: Theorem 1 provides a bound on the size of the
swarm, but the bound is conservative. This is because we enlarged
$\dot{V}$, and we used $e^{iT}(x^i-x^j)\leq \|x^i-x^j\|\|e^i\|$
and also assumed that the functions
$\exp\bigg(-\frac{\|x^i-x^j\|^2}{c}\bigg)\|x^i-x^j\|$ were at
their maximum value for all $i$ and $j$. Therefore, the actual
size of the swarm is, in general, much smaller than $\Omega$.

\textbf{Remark 5}: Under the assumption of
$\sum_{j=1}^Nw_{ij}=\sum_{j=1}^Nw_{ji}$, we obtain that the motion
of the swarm center only depends on the repulsion between the
swarm members.

\section{Extensions}

In Sections 2 and 3 we consider a specific function $f(y)$ as
defined in (2). In this section, we will consider a more general
function $f(y)$ that satisfies some assumptions. $f(y)$ is still
the social potential function that governs the interindividual
interactions and is assumed to have a long rang attraction and
short range repulsion nature. Following \cite{V. Gazi 1}, we make
the assumptions on the social potential function:

\textbf{Assumption 1}. The attraction/repulsion function
$f(\cdot)$ is of the form
\begin{equation}
f(y)=-y[f_a(\|y\|)-f_r(\|y\|)], y\in{R^n},
\label{eq7}
\end{equation}
where $f_a:R_+\rightarrow R_+$ represents (the magnitude of)
attraction term and has a long range, whereas $f_r:R_+\rightarrow
R_+$ represents (the magnitude of) repulsion term and has a short
range, and $R_+$ stands for the set of nonnegative real numbers,
$\|y\|=\sqrt{y^Ty}$ is the Euclidean norm.

\textbf{Assumption 2}. There are positive constants $a, b$ such
that for any $y\in R^n$,
\begin{equation}
f_a(\|y\|)=a, \ \ \ \  f_r(\|y\|)\leq \frac{b}{\|y\|}. \label{eq8}
\end{equation}

That is, we assume a fixed linear attraction function and a
bounded repulsion function.

\textbf{Theorem 2}: Consider the swarm described by the model in
(1) with an attraction/replusion function $f(\cdot)$ as given in
(7) satisfied (8). Then, all agents will eventually enter into and
remain in the bounded region
\begin{equation}
\Omega^*=\bigg\{x:
\sum_{i=1}^N\|x^i-\overline{x}\|^2\leq{\rho^2}\bigg\}, \label{eq9}
\end{equation}
where $\rho=\frac{4bM}{a\lambda_2};$ and $\lambda_2$ and $M$ are
defined as in Theorem 1; $\Omega^*$ provides a bound on the
maximum ultimate swarm size.

Following the proof of Theorem 1, Theorem 2 can be proved
analogously.

\section{Simulations}

In this section we will present some numerical simulations for the
nonreciprocal swarm described by Eqs. (1) and (2) in order to
illustrate the theory obtained in the previous section.

In these simulations we used the $f(\cdot)$ function which is
taken in the form of Eq. (2) with $a=1, b=20$, and $c=0.2$. The
coupling matrix $W$ is generated randomly and satisfies the before
conditions and assumptions.

Figs. 1-2 and Figs. 5-6 separately show the trajectories of the
swarm members and the swarm center in which there are $N=10$
individuals, and the four simulations run for 30s. In order to
more clearly describe the motion of the swarm members and the
swarm center, we also present the simulations of a five-agent
swarm, that is, Figs. 3-4, and the two simulations run for 100s.
It can be seen from Figs. 1-6 that at the beginning phase of the
simulations of the swarm member trajectories, all of the members
gradually aggregate and form a cohesive cluster. Then, they
continuously move in the same direction as a group, and eventually
evolve into an expending spiral motion as time increases.

\section{Conclusions}

In this paper, we have considered an anisotropic swarm model and
analyzed its aggregation. The model given here is a generalization
of the model in \cite{V. Gazi and K. M. Passino} and \cite{T.
Chu}. And the model is more applicable to the reality and is more
meaningful.

\end{multicols}
\end{document}